\documentclass[draftcls, onecolumn]{IEEEtran}
\usepackage{verbatim}
\usepackage{amssymb}
\usepackage[cmex10]{amsmath}
\usepackage{graphicx}
\usepackage{epsfig}
\usepackage{cite}
\usepackage{psfrag}
\usepackage{comment}
\usepackage{bm}
\usepackage{bbm}
\usepackage{url}
\usepackage{subfigure}
\usepackage{color}

\newcommand{\R}{\mathbb{R}}

\newcommand{\al}{\alpha}

\newcommand{\bd}{\mathbf}

\newcommand{\tl}{\tilde}

\newcommand{\C}{\mathbb{C}}

\newcommand{\sdp}{\succcurlyeq}
\newcommand{\pd}{\succ}
\newcommand{\bma}{\begin{bmatrix}}
\newcommand{\ebma}{\end{bmatrix}}

\newcommand{\nn}{\nonumber}

\newcommand{\mc}{\mathcal}

\newcommand{\vbar}{\overline{\bd v}}

\DeclareMathOperator{\diag}{diag}
\DeclareMathOperator{\diago}{diag}
\DeclareMathOperator{\Real}{Re}
\DeclareMathOperator{\Imag}{Im}

\DeclareMathOperator{\conj}{conj}

\newtheorem{defn}{Definition}
\newtheorem{thm}{Theorem}
\newtheorem{claim}[thm]{Claim}
\newtheorem{lem}[thm]{Lemma}

\begin{document}

\title{A Local Control Approach to Voltage Regulation in Distribution Networks}

\author{{Baosen Zhang, Alejandro D. Dom\'{i}guez-Garc\'{i}a and David Tse}
\thanks{B. Zhang and D. Tse are with the department of Electrical Engineering and Computer Sciences, University of California at Berkeley, Berkeley, CA, 94074 emails:\{zhangbao,dtse\}@eecs.berkeley.edu}
\thanks{A. D. Dom\'{i}guez-Garc\'{i}a is with the department of Electrical and Computer Engineering at University of Illinois at Urbana-Champaign, Urbana, IL 61801 email:aledan@illinois.edu}%
} 

\maketitle

\begin{abstract}
 This paper address the problem of voltage regulation in power distribution networks with deep penetration of distributed energy resources (DERs) without any explicit communication between the buses in the network. We cast the problem as an optimization problem with the objective of minimizing the distance between the bus voltage magnitudes and some reference voltage profile. We present an iterative algorithm where each bus  updates the reactive power injection provided by their DER. The update at a bus only depends on the voltage magnitude at that bus, and for this reason, we call the algorithm a local control algorithm. We provide sufficient conditions that guarantee the convergence of the algorithm and these conditions can be checked a priori for a set of feasible power injections. We also provide necessary conditions establishing that longer and more heavily loaded networks are inherently more difficult to control. We illustrate the operation of the algorithm through case studies involving 8-,34- and 123-bus test distribution systems. 
\end{abstract}

%
%

\section{Introduction}
Traditionally, power distribution systems have been engineered as radial networks that carry power from the utility feeder to the end consumers of electricity. These consumers are seen as a set of passive loads, and generally, the utility has no control of their electricity demand, and does not collect any real-time information from them. In the future, under the \emph{SmartGrid} initiative, distribution systems are envisioned to undergo dramatic transformation in structure and functionality. Instead of passive demands, buses in the network could include clusters of households with rooftop solar installations, electric vehicles, and storage; these generation and storage devices are commonly referred to as distributed energy resources(DERs). It has been shown that massive penetration of rooftop solar can cause voltages to rise in the distribution system, while large-scale deployment of electrical vehicles can cause voltages to drop \cite{Carvalho08,Guille09}. The objective of this paper is to address this voltage regulation problem. 


Voltage regulation in distribution networks is traditionally achieved via tap-changing under-load transformers and switched capacitors; however, while ensuring that voltage requirements are met, the setting of these devices might not optimal, in the sense that they may not minimize the thermal loss or some other objectives. On the other hand, by relying on a cyber-infrastructure for sensing, communication, and control, it is also possible to utilize DERs with capability to provide reactive power for optimal voltage regulation (see, e.g., \cite{Lam12} and the references therein). Unfortunately, the current speed of installation of DERs far outpaces the deployment of the cyber-infrastructure.\footnote{The advanced metering infrastructure (e.g. smart meters) transmit their information on a daily bases and are not used for real-time communication.} For example, most distribution networks in California have not implemented real-time two-way communication between DERs and the utilities, but the state already has 2900 MW of installed capacity \cite{SEIA13}, enough for about 600,000 households. Therefore we are faced with a situation where DERs are already in the distribution network with no communication network to coordinate them. Thus, an urgent question arises: in order to utilize DERs for voltage regulation in distribution networks, how should they be managed with no communication or very limited communication between the buses in the network?  

In our setting, we assume each bus can measure its own voltage magnitude, active and reactive power injections, but cannot communicate with other buses in the network. The available control action at a bus is the adjustment of its active and reactive power injections. Altering the active power injection at certain buses corresponds to demand response (DR) actions. Because of economical and regulatory considerations, dynamic DR is difficult to implement without a communication network; therefore, in this paper we only consider reactive power injections at each bus as the main control knob for DER-based voltage regulation. We assume that reactive power can be provided by, e.g., the power electronics grid interfaces of rooftop solar installations \cite{JoOo:00,LoKr:01,solar_bridge}. 

 In this paper, we propose a method for voltage regulation that relies on the local control actions\footnote{We call the control actions \emph{local control actions} to emphasize that only local information is available at each bus.} of reactive-power-capable DERs. We provide sufficient conditions such that local control is guaranteed to be successful in maintaining any specified voltage profile. We also provide conditions on the network topology and active power injection regions that shows when local control cannot maintain voltage stability in the system. We arrive at our results by casting the voltage regulation problem as an optimization problem, and investigate the conditions under which optimality can be achieved without communication between the buses. 

The problem of voltage regulation in distribution networks has received considerable attention in recent years. Centralized algorithms have been proposed in \cite{Baran07,Turitsyn10,Villacci06} and decentralized algorithms have been proposed in \cite{Farivar12,Lam12,Bolognani12,Robbins12}. For the latter, communication is generally assumed to be possible between subsets of the buses (e.g., \cite{Lam12} provides a distributed optimization algorithm that adheres to neighbor-to-neighbor communications, and \cite{Bolognani12} can choose different communication topology using an approximate model for power flow). In contrast, our proposed algorithm does not rely on a communication infrastructure. A related questions in context of voltage stability has been studied in \cite{Ilic:1986}. 

 This paper is organized as follows.  Section \ref{sec:model} introduces the notations and the power system distribution model adopted throughout the paper. Section \ref{sec:result} gives the main algorithm. Section \ref{sec:sr} provides convergence analysis of the algorithm. Section \ref{sec:sim} provides case studies,  and Section \ref{sec:con} concludes the paper. The proofs of the main results are provided in the Appendix. 

\section{Model and Problem Formulation} \label{sec:model}
This section introduces the power distribution system models and also provides the formulation of the problem. Consider a distribution network with $n$ buses. Today, most power distribution systems operate as radial networks connected to a single feeder\cite{Kersting06}; therefore we model the network as a connected tree rooted at the feeder. By convention, the feeder is labelled as bus $0$. The distance from bus $i$ to the root (bus 0) is the number of vertices on the unique path from bus $i$ to the root. The depth of the tree is the maximum of the bus distances to the root.    


Let $V_i=|V_i| e (j \theta_i)$ denote bus $i$ complex voltage, and define  vector $\bd v =[ V_0 \; V_1 \dots V_{n-1}]^T \in \C^n$. Let $I_i$, $P_i$ and $Q_i$ denote, the complex current, active and reactive power injections at bus $i$, respectively; and define $\bd i=[I_0 \; \dots \; I_{n-1}] \in \C^n$, $\bd p=[P_0 \; \dots \; P_{n-1}] \in \R^n$ ,and $\bd q=[Q_0 \; \dots \; Q_{n-1}] \in \R^n$. Let $\bd s=\bd p + j \bd q$ denote the complex power. The active powers in a distribution network are typically constrained by upper and lower capacity limits, i.e.,  of $\underline{P}_i \leq P_i \leq \overline{P}_i$, where the upper and lower bound are determined by by the types of devices at bus $i$. Let $\mc P =\{\bd p \in \R^{n-1}: \underline{P}_i \leq P_i \leq \overline{P}_i \}$ be the feasible region formed by the bus injections (see \cite{Zhang13,LTZ12} for more detailed discussions).

 Given two connected buses $i$ and $k$, let $y_{ik}=g_{ik}-jb_{ik}$ denote the admittance of the transmission line between them, and let $\bd Y$ be the bus admittance matrix \cite{Overbye04}. Then, the complex power injections are related to the voltage by $\bd s = \diag (\bd v \bd v^H \bd Y^H)$,
 where $(\cdot)^H$ denotes the Hermitian transpose, and the $\diag(\cdot)$ operator returns the diagonal if the argument is a square matrix or a diagonal matrix if the argument is a vector. We are interested in regulating the voltage in the distribution network, so let $V_0^{ref}, V_1^{ref},\dots,V_{n-1}^{ref}$ be the set of reference voltage magnitudes that we are interested in tracking. We assume that the feeder behaves as an infinite bus, thus $V_0=V_0^{ref}$. We assume that bus $i$ can control its reactive power, whereas the active power is set at a prescribed value $P_i^*$. We formulate the voltage regulation problem as the following feasibility problem:
 \begin{subequations} \label{eqn:feasible}
 \begin{align}
 \mbox{find } & \bd q \\
 \mbox{s.t. } & |V_i|=V_i^{ref},~ i=0,\dots,n-1 \\
 &  P_i=P_i^*, \; i=1,\dots,n-1 \label{eqn:p_f} \\
  & \bd p + j \bd q =\diag(\bd v \bd v^H \bd Y^H), \label{eqn:s_f}
\end{align}
\end{subequations}
 where \eqref{eqn:p_f} corresponds to the fact that active powers are fixed, and the equality constraints in (1d) correspond to the power flow equations.
That is, we seek to find a reactive power injection vector $\bd q$ such that the voltages are at their desired levels. The feasibility of \eqref{eqn:feasible} can be decided by the following optimization problem:
 \begin{subequations} \label{eqn:main_v} 
 \begin{align}
 \min_{\bd q }\; & \sum_{i=1}^{n-1} (|V_i|-V_i^{ref})^2 \label{eqn:obj_v} \\
 \mbox{s.t. } & P_i=P_i^*, \; i=1,\dots,n-1 \label{eqn:p_v} \\
 & V_0=V_0^{ref} \label{eqn:v0_v} \\
 & \bd p + j \bd q =\diag(\bd v \bd v^H \bd Y^H), \label{eqn:s_v}
 \end{align}
 \end{subequations}
 where the feeder bus (bus 0) always holds its voltage to be at $V_0^{ref}$. The objective function means that the system is trying to minimize the deviation of the voltage magnitudes to their reference values. The feasibility problem in \eqref{eqn:feasible} has a solution if and only if the optimal value of the optimization problem \eqref{eqn:main_v} is 0. 
 
 It turns out the objective function in \eqref{eqn:obj_v} is not easy to analyze. We rewrite \eqref{eqn:main_v} with a new objective function as:
  \begin{subequations} \label{eqn:main} 
  \begin{align}
  \min_{\bd q }\; & \sum_{i=1}^{n-1} \left(|V_i|^2-(V_i^{ref})^2 \right)^2 \label{eqn:obj} \\
  \mbox{s.t. } & P_i=P_i^*, \; i=1,\dots,n-1 \label{eqn:p} \\
  & V_0=V_0^{ref} \label{eqn:v0} \\
  & \bd p + j \bd q =\diag(\bd v \bd v^H \bd Y^H),  \label{eqn:s}
  \end{align}
  \end{subequations} 
 Intuitively, the reason that \eqref{eqn:obj} is much easier to handle than \eqref{eqn:obj_v} is that $|V_i|$ is not differentiable with respect to $\bd q$ but $|V_i|^2$ is. The rest of the paper presents a local control algorithm to solve \eqref{eqn:main}. 
  
\textbf{Remark:} In this paper we do not consider the constraints on bus reactive power injections $\bd q$. In an earlier work, where limited communication is present, all constraints were considered \cite{Lam12}. The next step of this work is to incorporate the constraints on reactive power injections. \hfill $\Box$ 

\subsection{Uniqueness of Power Flows} \label{sec:2A}
Given a power system with $n$ buses, the state of the system is described by $2n-1$ real numbers: the $n$ voltage magnitude and $n-1$ angles (taking bus $0$ to be the reference bus). This suggests that there are $2n-1$ degree of freedom in the power system. However, for transmission networks, specifying the voltage magnitudes and active powers $P_1,\dots,P_{n-1}$ is not enough to determine the state of the system since there could be multiple angles that result in those active powers \cite{BeVi:00}. In contrast, due to the tree structure of a distribution network \cite{LTZ12}, specifying voltage magnitudes and active powers at $n-1$ buses is enough to fully specify the state of the system. Therefore we often give the power system state in terms of voltage magnitude and active power; the corresponding angles and reactive powers are uniquely specified. 
  
\section{Local Control Algorithm} \label{sec:result}
Consider a distribution system that is operating at its reference voltage, i.e., with $|V_i|=V_i^{ref}$ for $i=0 \dots n-1$. Denote the corresponding active and reactive injection vectors by $\overline{ \bd p}$ and $\overline{\bd q}$, respectively. Due to changes in renewables and loads, the active power injections in the system may change from $\overline{\bd p}$ to $\bd p'$. If the reactive powers stayed constant at $\overline{\bd q}$, the voltage magnitudes are no longer at the reference values. This section proposes an algorithm to bring the voltage magnitudes back to the reference values. By adjusting the reactive power, the algorithm relies  purely on local voltage measurements and does not involve communication among nodes in the network. 

To correct the voltage magnitudes, there are $n-1$ control variables: the reactive power injections at all the buses except the feeder (bus $0$). We assume that the problem in \eqref{eqn:feasible} has a solution, or equivalently the optimal value of the problem in \eqref{eqn:main} is 0, and we design a local control algorithm that finds the optimal reactive power injection vector $\bd q$. Since bus $i$ can only measure local quantities, we propose the following iterative update algorithm:

\textbf{Local Control Algorithm:}
\begin{itemize}
\item Initialize the system at $t=0$. 
\item At time $t$, bus $i$ measures $|V_i[t]|$ and makes the following update: 
\begin{equation} \label{eqn:update_i}
Q_i[t+1]=Q_i[t]-d_i \left(|V_i[t]|^2-(V_i^{ref})^2 \right)
\end{equation}
for $i=1,\dots,n-1$. Intuitively, a bus should inject more reactive power to raise its voltage, so $d_i$ is assume to be positive (see \cite{Robbins12} for a deeper discussion). 
\item $t \leftarrow t+1$ and repeat
\end{itemize}
By a slight abuse of notation, let $\bd q=[Q_1 \; \dots Q_{n-1}]^T$ since we do not explicitly control the reactive injection at the feeder. The matrix version of the update in \eqref{eqn:update_i} can be written as
\begin{equation} \label{eqn:update}
\bd q[t+1]= \bd q[t]- \bd D \left(|\bd v[t]|^2 - (\bd v^{ref})^2 \right),
\end{equation}
where $\bd D=\diag(d_1,\dots,d_{n-1})$ is a diagonal matrix and the squares are taken component-wise. Note both $|\bd v[t]|$ and $\bd v^{ref}$ are both vectors. We envision the local control algorithm would run in a continuous fashion to supplement the actions of devices such as tap-changing under-load transformers (TCULs) and switched capacitors. This local control algorithm is similar to the first-stage algorithm in \cite{Robbins12} where we replaced $|V_i|$ by $|V_i|^2$.
 
\section{Convergence of the Local Control Algorithm}
The convergence properties of the local control algorithm are investigated in this section. Since the active bus injections are the parameters in the voltage regulation problem, the convergence of the algorithm also depends on the prescribed active power injections. Recall from Section \ref{sec:model} that $\mc P$ is the feasible injection region determined by the active bus power constraints. In this section, we answer the following two questions:
\begin{enumerate}
\item Under what condition can we ensure that the algorithm converges for every $\bd p \in \mc P$?
\item How does network topology influence the convergence of the algorithm? That is, for the two networks in Fig. \ref{fig:star_line} and same active injection vector, does the algorithm converge for one of the networks but not the other one? 
\begin{figure}[ht]
\centering
\subfigure[Star Topology]{
\includegraphics[scale=0.7]{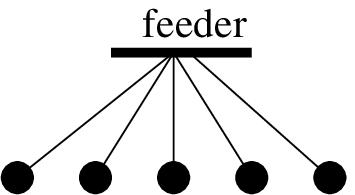}}
\subfigure[Line Topology]{
\includegraphics[scale=0.7]{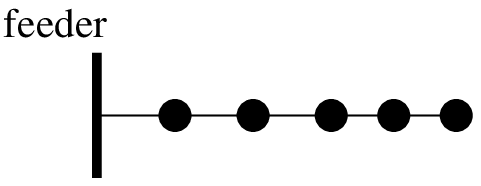}}
\caption{We are interested in how the topology of the network affect the convergence property of the local control algorithm. We show that it is easier for the local algorithm to converge on the star network than the tree network. Therefore the depth of the network limits the convergence of the algorithm.}
\label{fig:star_line}
\vspace{-0cm}
\end{figure}
\end{enumerate}

Since power system planners generally plan for the worst-case scenario, answering the first question would ensure that the algorithm would work for every possible $\bd p$ and communication is not needed to regulate the voltage. Theorem \ref{thm:sr} below states that it is in fact sufficient to check the convergence of the algorithm at a single $\bd p$ to ensure the converge for the entire feasible injection region $\mathcal{P}$.

In section IV-C, we show that the local control algorithm converges for the star topology under a much broader set of conditions than for the line topology. Since typical distribution networks are more like lines than stars, they are in some sense inherently more difficult to control only using local actions. Before stating the convergence theorems for the entire $\mc P$, we first state the result for a single power injection vector. 
\subsection{Convergence for a Power Injection Vector }
The local control algorithm shows that we can view the $n-1$ voltages on the non-feeder buses as a function of the $n-1$ reactive power injections at these buses. The active powers and the feeder voltage can be interpreted as parameters in this function. From this point of view, we make the following definition:
Let $\bd M$ be the $n-1$ by $n-1$ Jacobian matrix where the $(i,k)$-th entry is
\begin{equation} \label{eqn:M}
M_{ik}= \frac{\partial |V_i|^2}{\partial Q_k}; 
\end{equation}
the partial derivative is evaluated at some active power and voltage magnitudes. 

Note that $\bd M$ depends on both values of both the active powers and voltage magnitudes; however,  the variation of $\bd M$ with respect to voltage magnitude is relatively small over a wide range of operating conditions \cite{BeVi:00}. In this regard, we will always evaluate $\bd M$ at some $\bd p$ and the reference voltage magnitudes. To emphasize the dependence of $\bd M$ on $\bd p$, we write it as $\bd M (\bd p)$.  

The next lemma gives the sufficient condition when the local control algorithm converges. 
\begin{lem} \label{lem:1}
Suppose the optimization problem in \eqref{eqn:main} has optimal value $0$ (that is, the problem in \eqref{eqn:feasible} is feasible). Let $\bd p$ be the active power injection vector. The update in \eqref{eqn:update_i} converges if $\bd D \bd M(\bd p) + \bd M(\bd p)^T \bd D \pd 0 $, where $ \pd 0 $ denotes positive definiteness.
\end{lem}
The proof is an application of Prop. 2.1 in \cite{Bertsekas97} and is given in Appendix \ref{app:1}. 

\subsection{Stability Region} \label{sec:sr}
The matrix inequality $\bd D \bd M(\bd p) + \bd M(\bd p)^T \bd D \pd 0 $ is in the form of a Lyapunov equation, and this motivates the next definition 
\begin{defn}
An active power injection vector is called stable if there exists a diagonal matrix $\bd D$ such that $\bd D \bd M(\bd p) + \bd M(\bd p)^T \bd D \pd 0$. The set of all stable active power injections is called the stability region and is denoted by $\mc S_p$. 
\end{defn} 

Deciding the existence of a stabilizing $\bd D$ for a given $\bd p$ is a convex problem, and can be solved using semidefinite programming techniques \cite{Boyd04}. Theorem \ref{thm:sr} below shows how the stability for the entire region $\mc P$ can be checked. 

Given the feasible injection region $\mc {P}$, let $\bd p_{\min}$ be the minimum injection vector in $\mc P$. That is, $\bd p_{\min}=[\underline{P}_1 \; \dots \; \underline{P}_{n-1}]$ (recall $\underline{P_i}<0$ is the lower bound on the active power injection, which correspond to the maximum demand). The next theorem gives a sufficient condition for the stability of $\mc P$ in terms of the stability of $\bd p_{\min}$. 
\begin{thm} \label{thm:sr}
Let $\bd p_{\min}$ be the minimum injection vector.  Let $\bm \theta_{\min}$ be the corresponding angles when the voltage magnitudes are at their reference values (recall Section \ref{sec:2A}). Suppose $|\theta_{\min,ik}| \leq \tan^{-1}(b_{ik}/g_{ik})$. Then if $\bd M(\bd p_{\min})$ is stable, $\bd M(\bd p)$ is stable for all $\bd p \in \mc P$. That is, $\mc P \subseteq \mc S_p$.  
\end{thm}

The proof of Theorem \ref{thm:sr} is given in the Appendix. The significance of the theorem is that it suffices to check the stability of a single point, and this would ensure the global stability of the entire $\mc P$. Thus, in practice it is easy to check whether $\mc P \subseteq \mc S_p$ and how much demand need to be reduced to achieve stability under local controls. 

The condition on $\theta_{\min,ik}$ is discussed in detail in \cite{LTZ12}, and is equivalent to the thermal constraints on the transmission lines. They are expected to hold in almost all practical situations. 

\subsection{Effects of Topology on the Stability Region} \label{sec:long}
Theorem \ref{thm:sr} states that the stability of $\mc P$ is determined by the minimum injection vector. In this section, we show that the stability region $\mc S_p$ shrinks as the depth of the network increases. This shows that long networks are inherently more difficult to control, or equivalently, a line network is inherently more unstable than a star network. Due to space constraints, rather than stating a general theorem, we state a result pertaining to line networks. 
\begin{lem} \label{lem:long}
Consider a $n$-bus homogeneous line network where every line has admittance $1-j$. Then the stability region of $\mc S_p$ approaches the single point $\bd 0$ as $n$ goes to infinity. 
\end{lem}

This lemma states that as the depth of the network increases, communication becomes critical in maintaining the voltage stability of the network.  Practical networks often have a large depth, e.g., the 123-bus system has a depth of 23; therefore, some communication would likely be necessary for large distribution networks. 

\section{Case Studies} \label{sec:sim}
We analyze the stability of the proposed local control algorithm for three of the IEEE test feeders: the 8-bus, 34-bus and the 123-bus network (see \cite{testfeeders} for system data). Since only one demand point is given in the data, we compute $\bd M$ at that point. By Theorem \ref{thm:sr}, all demand vectors less than the given demand is stable. Table \ref{tab:stable} shows whether the networks is stable under its current demand.   
\begin{table}[ht]
\centering
\begin{tabular}{|c|c|c|c|}
Networks & 8-bus & 34-bus & 123-bus \\
\hline
Depth & 6 & 19 & 23 \\
Stable & yes & yes & no 
\end{tabular}  
\caption{Stability of the three network at the demand given in \cite{testfeeders}.}
\label{tab:stable}
\vspace{-0.5cm}
\end{table}

\begin{figure}[!t]
\centering
\includegraphics[scale=0.35]{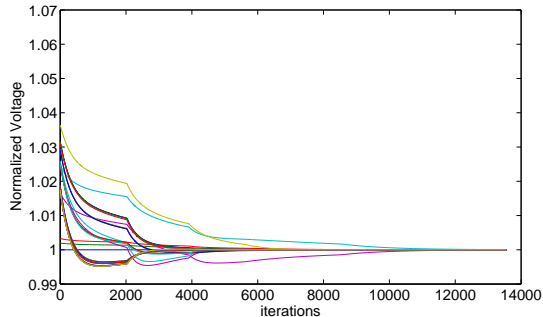}
\caption{The behavior of voltage magnitudes under the local control algorithm. Each curve is the normalized voltage, i.e., $|V_i[t]|/V_i^{ref}$. All of the curves converges to $1$, which means that all voltages are at their reference.}
\label{fig:34sim_a}
\end{figure}

\begin{figure}[!t]
\centering
\includegraphics[scale=0.35]{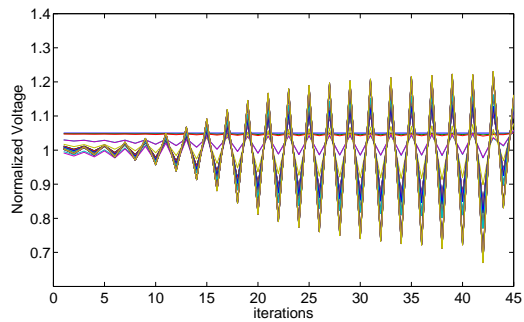}
\caption{The behavior of voltage magnitudes under the local control algorithm. Each curve is the normalized voltage, i.e., $|V_i[t]|/V_i^{ref}$. The demand in the system is 3.5 times the demand in Fig. \ref{fig:34sim_a}. Due the heavier loading, the network becomes unstable under the local control algorithm and experiences large voltage swings.}
\label{fig:34sim_b}
\vspace{-0.7cm}
\end{figure}

We focus on the 34-bus network and demonstrate the action of the local control algorithm. Figures \ref{fig:34sim_a} and \ref{fig:34sim_b} respectively show the voltage magnitudes in the network when the local action is able to stabilize the voltage and when it fails to do so. In Fig. \ref{fig:34sim_a} the demand is set to be the current demand in \cite{testfeeders}, and initially, we vary the voltage randomly by 10 \% and apply the local control algorithm in Section \ref{sec:result}. The voltage magnitudes converge to the reference voltages.

The local control algorithm fails to converge once the demand in the network is increased to 3 times its current value. In Fig. \ref{fig:34sim_b}, we increase the demand to $3.5$ times the current demand, and again vary the voltage randomly by 10 \% from its reference value; the figure shows the unstable behavior of the voltage magnitudes.

\section{Conclusion} \label{sec:con}
In this paper we investigated the utilization of reactive-power-capable DERs for voltage regulation in distribution networks without explicit communications among the buses. We proposed an algorithm where each bus controls its reactive power injections sensing only its own voltage magnitude. We provided sufficient conditions when the algorithm converges for the entire feasible active power injection region of a network, and necessary conditions that shows longer and heavier loaded networks are harder to control using localized algorithms. We validated our claims in case studies using the IEEE 34-bus test feeder. 

\bibliographystyle{IEEEtran}
\bibliography{mybib}

\begin{thebibliography}{10}
\providecommand{\url}[1]{#1}
\csname url@samestyle\endcsname
\providecommand{\newblock}{\relax}
\providecommand{\bibinfo}[2]{#2}
\providecommand{\BIBentrySTDinterwordspacing}{\spaceskip=0pt\relax}
\providecommand{\BIBentryALTinterwordstretchfactor}{4}
\providecommand{\BIBentryALTinterwordspacing}{\spaceskip=\fontdimen2\font plus
\BIBentryALTinterwordstretchfactor\fontdimen3\font minus
  \fontdimen4\font\relax}
\providecommand{\BIBforeignlanguage}[2]{{%
\expandafter\ifx\csname l@#1\endcsname\relax
\typeout{** WARNING: IEEEtran.bst: No hyphenation pattern has been}%
\typeout{** loaded for the language `#1'. Using the pattern for}%
\typeout{** the default language instead.}%
\else
\language=\csname l@#1\endcsname
\fi
#2}}
\providecommand{\BIBdecl}{\relax}
\BIBdecl

\bibitem{Carvalho08}
P.~Carvalho, P.~Correia, and L.~Ferreira, ``Distributed reactive power
  generation control for voltage rise mitigation in distribution networks,''
  \emph{{IEEE} Trans. Power Sys.}, vol.~23, no.~2, pp. 766--772, May 2008.

\bibitem{Guille09}
C.~Guille and G.~Gross, ``A conceptual framework for the vehicle-to-grid (v2g)
  implementation,'' \emph{Energy Policy}, vol.~37, no.~11, pp. 4379--4390, Nov.
  2009.

\bibitem{Lam12}
\BIBentryALTinterwordspacing
A.~Y. Lam, B.~Zhang, A.~Dominguez-Garcia, and D.~Tse. (2012) Optimal
  distributed voltage regulation in power distribution networks. [Online].
  Available: \url{http://arxiv.org/abs/1204.5226}
\BIBentrySTDinterwordspacing

\bibitem{SEIA13}
\BIBentryALTinterwordspacing
{Solar Energy Industry Association}. (2013) Facts on the california solar
  industry. [Online]. Available:
  \url{http://www.seia.org/state-solar-policy/california}
\BIBentrySTDinterwordspacing

\bibitem{JoOo:00}
G.~Joos, B.~Ooi, D.~McGillis, F.~Galiana, and R.~Marceau, ``The potential of
  distributed generation to provide ancillary services,'' in \emph{Proc. of
  IEEE Power Engineering Society Summer Meeting}, Seattle, WA, 2000.

\bibitem{LoKr:01}
D.~Logue and P.~Krein, ``Utility distributed reactive power control using
  correlation techniques,'' in \emph{Proc. of IEEE Applied Power Electronics
  Conference}, Anaheim, CA, March 2001.

\bibitem{solar_bridge}
\BIBentryALTinterwordspacing
\text{SolarBridge Technologies}. (2009) \text{Pantheon Microinverter}. South
  Plainfield, NJ. [Online]. Available: \url{http://www.petrasolar.com/}
\BIBentrySTDinterwordspacing

\bibitem{Baran07}
M.~Baran and I.~El-Markabi, ``A multiagent-based dispatching scheme for
  distributed generators for voltage support on distribution feeders,''
  \emph{IEEE Trans. Power Syst.}, vol.~22, no.~1, pp. 52--59, Feb. 2007.

\bibitem{Turitsyn10}
K.~Turitsyn, P.~Sulc, S.~Backhaus, and M.~Chertkov, ``Distributed control of
  reactive power flow in a radial distribution circuit with high photovoltaic
  penetration,'' in \emph{Proc. of {IEEE} Power and Energy Society General
  Meeting, 2010}, Minneapolis, MN, July 2010, pp. 1--6.

\bibitem{Villacci06}
D.~Villacci, G.~Bontempi, and A.~Vaccaro, ``An adaptive local learning-based
  methodology for voltage regulation in distribution networks with dispersed
  generation,'' \emph{{IEEE} Trans. Power Sys.}, vol.~21, no.~3, pp.
  1131--1140, Aug. 2006.

\bibitem{Farivar12}
M.~Farivar, R.~Neal, C.~Clarke, and S.~Low, ``Optimal inverter var control in
  distribution systems with high pv penetration,'' in \emph{Proc. of the IEEE
  Power and Engergy Society General Meeting}, San Diego, CA, June 2012.

\bibitem{Bolognani12}
S.~Volognani and S.~Zampieri, ``A distributed control strategy for reactive
  power compensation in smart microgrids,'' \emph{{ArXiv}:1106.5626v3
  [math.OC]}, 2012.

\bibitem{Robbins12}
B.~A. Robbins, C.~N. Hadjicostis, and A.~D. Dom\'{i}nguez-Garc\'{i}a, ``A
  two-stage distributed architecture for voltage control in power distribution
  systems,'' \emph{to Appear in IEEE Transactions on Power Systems}, 2012.

\bibitem{Ilic:1986}
M.~Ilic-Spong, J.~Thorp, and M.~Spong, ``Localized response performance of the
  decoupled \text{Q-V} network,'' \emph{{IEEE} Trans. Circuits and Systems},
  vol.~33, no.~3, pp. 316 -- 322, March 1986.

\bibitem{Kersting06}
W.~H. Kersting, \emph{Distribution system modeling and analysis}.\hskip 1em
  plus 0.5em minus 0.4em\relax CRC Press, 2006.

\bibitem{Zhang13}
B.~Zhang and D.~Tse, ``Geometry of injection regions of power networks,''
  \emph{IEEE Transactions on Power Systems}, vol.~28, no.~2, pp. 788--797,
  2013.

\bibitem{LTZ12}
J.~Lavaei, D.~Tse, and B.~Zhang, ``Geometry of power flows in tree networks,''
  in \emph{Proc. of IEEE PES General Meetings}, 2012.

\bibitem{Overbye04}
T.~J. Overbye, X.~Cheng, and Y.~Sun, ``A comparison of the {AC} and {DC} power
  flow models for lmp calculations,'' in \emph{Proceedings of the 37th Hawaii
  International Conference on System Sciences}, 2004.

\bibitem{BeVi:00}
A.~Bergen and V.~Vittal, \emph{Power System Analysis}.\hskip 1em plus 0.5em
  minus 0.4em\relax Upper Saddle River, NJ: Prentice Hall, 2000.

\bibitem{Bertsekas97}
D.~P. Bertsekas and J.~N. Tsitsiklis, \emph{Parallel and Distributed
  Computation: Numerical Methods}.\hskip 1em plus 0.5em minus 0.4em\relax
  Athena Scientific, 1997.

\bibitem{Boyd04}
S.~Boyd and L.~Vandenberghe, \emph{Convex Optimization}.\hskip 1em plus 0.5em
  minus 0.4em\relax Cambridge, 2004.

\bibitem{testfeeders}
{Distribution Test Feeder Working Group}, ``Distribution test feeders,''
  \\http://ewh.ieee.org/soc/pes/dsacom/testfeeders/index.html, 2010.

\end{thebibliography}
\vspace{-0.3cm}
\appendix[Proof of Lemma \ref{lem:1}] \label{app:1}
\begin{proof}
The proof is an application of the following lemma from \cite{Bertsekas97} (Prop. 2.1, Chap. 3):
\begin{lem} \label{lem:ber}
Consider the problem of 
\begin{equation*}
\min_x \; F(x).
\end{equation*}
Suppose $F$ is continuously differentiable and bounded from below, and there exists a constant $K$ such that 
\begin{equation*}
||\nabla F(x)-\nabla F(y)||_2 \leq K ||x-y||_2.
\end{equation*}
Then the updates
\begin{equation*}
x[t+1]=x[t]+\al c[t]
\end{equation*}
converges to $x^*$ such that $\nabla F(x^*)=0$ for small enough $\al$ if there exists positive constants $K_1$ and $K_2$ such that
\begin{subequations} 
\begin{align}
& ||c[t]||_2 \geq K_1 ||\nabla F(x[t])||_2 \label{eqn:K1}\\ 
& c[t]^T \nabla F(x[t]) \leq -K_2 ||c[t]||_2^2. \label{eqn:K2}
\end{align}
\end{subequations}
\end{lem}

To apply Lemma \ref{lem:ber}, note the objective function of \eqref{eqn:main} is a function of $\bd q$.  By standard argument from complex analysis, the objective function is  bounded from below, continuously differentiable and its gradient is Lipschitz ($K$ exists). Note $\nabla F(\bd q)=\bd M$.  Let 
\begin{equation}
\bd c[t]=-\bd D ( |\bd v[t]|^2-(\bd v^{ref})^2). 
\end{equation}

As long as $\bd M$ is not identically $0$, we can find a sufficiently small constant $K_1$ such that 
\begin{align*}
& ( |\bd v[t]|^2-(\bd v^{ref})^2)^T \bd D^T \bd D ( |\bd v[t]|^2-(\bd v^{ref})^2) \\
\geq & K_1 ( |\bd v[t]|^2-(\bd v^{ref})^2)^T \bd M^T \bd M ( |\bd v[t]|^2-(\bd v^{ref})^2),
\end{align*}
therefore satisfying \eqref{eqn:K1}. 

The condition in \eqref{eqn:K2} becomes: there exists $K_2>0$ such that
\begin{align*}
& ( |\bd v[t]|^2-(\bd v^{ref})^2)^T \bd D^T \bd M \\
\geq & K_2 ( |\bd v[t]|^2-(\bd v^{ref})^2)^T \bd D^T \bd D ( |\bd v[t]|^2-(\bd v^{ref})^2). 
\end{align*}
Rearranging, we need to find a constant $K_2>0$ such that $\bd D^T \bd M- K_2 \bd D^T \bd D$ is positive definite. Since $\bd D^T \bd D$ is positive definite\footnote{A square matrix $\bd A$, not necessarily symmetric, is said to be positive definite if $\bd x^T \bd A \bd x>0$ for all $\bd x \neq 0$. If $\bd A$ is symmetric, we write $\bd A \pd 0$.}, $K_2$ exists if $\bd D^T \bd M$ is positive definite, or $\bd D^T \bd M + \bd M^T \bd D \pd 0$. 
\end{proof}

\appendix[Proof of Theorem \ref{thm:sr} and Lemma \ref{lem:long}] \label{sec:proof}
\vspace{-0.1cm}
Central to the proofs of Theorem \ref{thm:sr} and Lemma \ref{lem:long} is an approximate analytical expression for $\bd M$. Following the notation in \cite{Bolognani12}, let $\bd e_1$ be the $n$ by 1 vector with $1$ at the first component and $0$'s everywhere else. Define the matrix $\tl{\bd X}$ to be a matrix such that 
\begin{align*}
\bd{\tl X} \bd Y &=\bd I - \bd 1 \bd e_1^T \\
\bd{\tl X} \bd e_1 &= \bd 0, 
\end{align*}
where $\bd 1$ is the all $1$'s vector. Lemma 1 from \cite{Bolognani12} states that $\bd {\tl X}$ is symmetrical and uniquely given by 
\begin{equation}
\bma \bd {\tl X} & \bd 1 \\ \bd 1^T & 0 \ebma = \bma \bd Y & \bd e_1 \\ \bd e_1 & 0 \ebma^{-1}. 
\end{equation}
By the fact that $\tl {\bd X}$ is symmetric and $\tl{\bd X} \bd e_1 =\bd 0$, the first row and first column of $\tl{\bd X}$ are all $0$, and we partition $\tl{\bd X}$ as 
\begin{equation}
\tl{\bd X} = \bma 0 & 0 \\ 0 & \bd X \ebma,
\end{equation}
where $\bd X$ is the lower $n-1$ by $n-1$ sub-matrix. 
Given vectors $\bd x$ and $\bd y$ of the same length, we write $\bd x / \bd y$ to mean component-wise division. Given a complex matrix $\bd A$, $\overline{\bd A}$ denote the component-wise conjugation. 

The following theorem approximates $\bd M$ in terms of $\bd X$. 
\begin{thm}\label{thm:M}
Given an $n$ bus network, let $\bd v$ be the set of complex voltages on buses $2$ to $n$. Let $q$ be the corresponding reactive powers; then 
\begin{align} \label{eqn:MX}
\bd M = & 2\{ \Imag [\diag (\overline{\bd v}) \bd X \diag (1/ \overline{\bd v})] \\&- \Real[\diag (\overline{\bd v}) \bd X \diag (\bd q / \overline{\bd v}^2) \bd X^H \diag (1/\bd v)]\} \nn \\
& +O(1/\bd v^3). \nn
\end{align}
\end{thm}
Since the bus voltage magnitudes in the distribution system are typically at least $2.4$~kV, the last term in \eqref{eqn:MX} is negligible and can be dropped; therefore, $\bd M$ is approximately
\begin{align} \label{eqn:M1}
\bd M \approx  & 2\{ \Imag [\diag (\overline{\bd v}) \bd X \diag (1/ \overline{\bd v})] \\ 
&-\Real[\diag (\overline{\bd v}) \bd X \diag (\bd q / \overline{\bd v}^2) \bd X^H \diag (1/\bd v)]\ \nn.
\end{align}
The proof of Theorem \ref{thm:M} is given in Section \ref{sec:proofM}. The rest of this section uses Theorem \ref{thm:M} to prove Theorem \ref{thm:sr} and Lemma \ref{lem:long}. 

\subsection{Proof of Theorem \ref{thm:sr}}
First we show: 
\begin{claim}
Given two vectors $\bd q  \leq \tl{\bd q}$. Suppose $\bd M(\tl{\bd q})$ is stable, then $\bd M (\bd q)$ is diagonally stable. 
\end{claim}
Let $\bd D$ be the stabilizing matrix of $\bd M(\tl{\bd q})$. By \eqref{eqn:M1}, 
\begin{align}
&1/2 (\bd D \bd M(\bd q)+ \bd M(\bd q)^T \bd D) \nn \\
 &= 1/2 (\bd D \bd M(\tl{\bd q}-\tl{\bd q}+\bd q) + \bd M(\tl{\bd q}-\tl{\bd q}+\bd q)^T \bd D) \nn \\
&= 1/2 (\bd D \bd M(\tl{\bd q})+ \bd M(\tl{\bd q})^T \bd D) \nn \\
&+ 1/2  \bd D  \Real \left( \diag(\vbar) \bd X \diag ((\tl{\bd q} -\bd q)/ \vbar^2) \bd X^H \diag(1/\bd v) \right ) \nn \\
&+1/2 \Real \left( \diag(\vbar) \bd X \diag ((\tl{\bd q} -\bd q)/ \vbar^2) \bd X^H \diag(1/\bd v) \right )^T \bd D \nn \\
&\stackrel{(a)}{\sdp} 0 \nn,  
\end{align} 
where $(a)$ follows from $\tl{\bd q} \geq \bd q$. 

The next claim finishes the proof. 
\begin{claim}
Let $\bd p$ ($\bd q$) and $\tl{\bd p}$ ($\tl{\bd q}$) be two bus active power (reactive power) injection vectors. Suppose the usual angle constraints are satisfied. If $\tl{\bd p} \geq \bd p$, then $\tl{\bd q} \leq \bd q$. 
\end{claim}
We illustrate the proof for a line network. 
Let $\bd \theta$ ($\tl{\bd \theta}$) be the angles corresponding to $\bd p$ ($\tl{\bd p}$). Consider the active power injection at bus $n-1$. Since $ 0> \tl{P}_{n-1} >P_{n-1}$ and $\theta_{n-1}$ satisfies the angle constraints, $\tl{\theta}_{n-1} < \theta_{n-1}$ and $\tl{Q}_{n-1} < Q_{n-1}$. Repeating the argument up the tree yields the claim. 
\vspace{-0.2cm}  
\subsection{Proof of Lemma \ref{lem:long}}
The $\bd X$ matrix is 
\begin{equation*}
\bd X = 
\frac{1+j}{2}
\bma
1 & 1 & 1 & \dots & 1 \\
1 & 2 & 2 & \dots & 2 \\
1 & 2 & 3 & \dots & 3 \\
  &   &   & \vdots &  \\
1 & 2 & 3 & \dots & n-1 
\ebma
\end{equation*}

Since the gains $d_i$ are positive in \eqref{eqn:update_i}, a necessary condition for $\bd M$ to be diagonally stable for $\bd M$ to have a positive diagonal. From \eqref{eqn:M1},
\begin{equation*}
\diag{\bd M}= 2\{ \Imag[\diag (\bd X)] - \Real[\diag(\bd X \bd X^H) \odot \diag(q/|\bd v|^2),
\end{equation*}
where $\odot$ denotes component-wise multiplication. Simple calculation gives $\Imag(X_{ii}) = \frac{i}{2}$
and
\begin{equation*}
\Real[(\bd X \bd X^H)_{ii}] = \frac{(i-1)i(2i-1)}{6} + (n-i)i^2.
\end{equation*}
Therefore the diagonal of $\bd X$ is positive only if 
\begin{equation}
\frac{q_i}{|V_i|^2} \leq \frac{3}{(i-1)(2i-1)+6(n-i)i}.
\end{equation}
Since all buses are withdrawing point, for at least one $i$, $q_i \geq 0$. So $q_i \rightarrow 0$ as $n \rightarrow \infty$. 
\vspace{-0.4cm}
\subsection{Proof of Theorem \ref{thm:M}} \label{sec:proofM}
The state of a power system can be represented by the following set of equations:
\begin{subequations}\label{eqn:Yv}
\begin{align}
\bd i = \bd Y \bd v \label{eqn:Yva}\\
\bd v \odot \overline{\bd i} = \bd s. \label{eqn:Yvb}
\end{align}
\end{subequations}
 Multiplying both sides of \eqref{eqn:Yva} by $\tl{\bd X}$ and looking at the voltages at bus $1$ to bus $n-1$, we have $ \bd v = \bd X \conj(\frac{\bd s}{\bd v}) + V_0 \bd 1$,
from which it follows
\begin{align}
\bd v & = \bd X ( \frac{\bd p -j \bd q}{\overline{\bd v}}) + V_0 \bd 1 \nn \\
&= -j\bd X (\bd q/\overline{\bd v})+  \bd X ( \bd p/\overline{\bd v}) + V_0 \bd 1 \nn
\end{align}

Let $\frac{\partial \bd v}{\partial \bd q}$ be the Jacobian of $\bd v$ with respect to $\bd q$, by linearity
\begin{equation}
\frac{\partial \bd v}{\partial \bd q}= -j \bd X \frac{\partial \bd q /\overline{\bd v}}{\partial \bd q}. 
\end{equation}
\begin{equation}
(\frac{\partial \bd q /\overline{\bd v}}{\partial \bd q})_{ik}=
\begin{cases}
 \frac{1}{\overline v_i} - \frac{q_i \partial \overline v_i / \partial q_i }{\overline v_i^2} & { if } \;\; i=k\\
 -\frac{q_i \partial \overline v_i / \partial q_k }{\overline v_i^2} & { if } \; \; i \neq k
\end{cases}
\end{equation}
which can be compactly written in matrix form as
\begin{equation}
\frac{\partial \bd q /\overline{\bd v}}{\partial \bd q}= \diag(1/\overline{v}_i) - \diag(\bd q / \overline{\bd v}_i^2) \overline{\frac{\partial \bd v}{\partial \bd q}}
\end{equation}
and
\begin{equation} \label{eqn:vq}
\frac{\partial \bd v}{\partial \bd q}=-j\bd X \left( \diag(1/\overline{v}) - \diag(\bd q / \overline{\bd v}^2) \overline{\frac{\partial \bd v}{\partial \bd q}} \right). 
\end{equation}
From \eqref{eqn:vq}, we can calculate $\frac{\partial \bd v}{\partial \bd q}$ by matrix inversion, but it is not useful in deriving an analytical formula for $\frac{\partial |\bd v|^2}{\partial \bd q}$; therefore, we do not explicit solve \eqref{eqn:vq}. Instead, we have
\begin{align}
\frac{\partial |\bd v|^2}{\partial \bd q} &= \frac{\partial \bd v \odot \overline{\bd v}}{\partial \bd q} \nn \\
&= \diag(\bd v) \overline{\frac{\partial \bd v}{\partial \bd q}} + \diag(\overline{\bd v})\frac{\partial \bd v}{\partial \bd q} \nn \\
&= \diag(\bd v) [j \overline{\bd X} (\diag(1/\bd v)-\diag(\bd q / \bd v^2) \frac{\partial \bd v}{\partial \bd q})] \nn \\
&+ \diag(\overline{\bd v}) [-j \bd X(\diag(1/\overline{\bd v})-\diag(\bd q / \overline{\bd v}^2) \overline{\frac{\partial \bd v}{\partial \bd q}})] \nn \\
&= 2 \Imag \left[ \diag(\overline{\bd v}) \bd X \diag(1/\overline{\bd v}) \right ] \nn \\
&-2 \Imag \left[\diago(\overline{\bd v})\bd X \diago(\bd q/ \overline{\bd v}^2)\overline{\frac{\partial \bd v}{\partial \bd q}} \right]. \label{eqn:2Im}
\end{align}
Substituting \eqref{eqn:vq} into the second term of \eqref{eqn:2Im}, we obtain
\begin{align}
& \Imag \left[\diago(\overline{\bd v})\bd X \diago(\bd q/ \overline{\bd v}^2)\overline{\frac{\partial \bd v}{\partial \bd q}} \right] \\
 =& \Imag \{ \diago(\overline{\bd v})\bd X \diago(\bd q/ \overline{\bd v}^2) \cdot \nn \\
 & \left(j \overline{\bd X} (\diag(1/\bd v) - \diag(\bd q/ \bd v^2) \frac{\partial \bd v}{\partial \bd q}) \right)   \} \nn \\
= & \Imag \left[j \diag(\overline{ \bd v}) \bd X \diag(\bd q/ \overline{\bd v}^2) \overline{\bd X} \diag(1/\bd v) \right] + O(1/\bd v^3) \nn \\
\stackrel{(a)}{=} & \Real \left[\diag(\overline{ \bd v}) \bd X \diag(\bd q/ \overline{\bd v}^2) \bd X^H \diag(1/\bd v) \right] + O(1/\bd v^3)  \label{eqn:Re}
\end{align}
where $(a)$ follows from the fact $\bd X$ is symmetric and $\Imag(jz)=\Real(z)$ for complex $z$.  
Combining \eqref{eqn:2Im} and \eqref{eqn:Re} gives Theorem~\ref{thm:M}.  Compared to \cite{Bolognani12}, \eqref{eqn:MX} is a more accurate version of $\bd M$. In \cite{Bolognani12}, $\bd M$ is simply approximated to be $2 \Imag[\bd X]$. 

\end{document}